\documentclass[reqno,12pt]{amsart}
\usepackage{latexsym,amsmath,amssymb,amsbsy,amscd,bm,amsfonts,mathrsfs}
\usepackage{graphicx}
\usepackage{hyperref}
\hypersetup{
  colorlinks   = true, 
  urlcolor     = blue, 
  linkcolor    = blue, 
  citecolor    = red   
}
\theoremstyle{plain}
\newtheorem{theorem}{Theorem}[section]

\newtheorem{lemma}{Lemma}[section]
\newtheorem{prop}{Proposition}[section]
\newtheorem{cor}{Corollary}[section]

\theoremstyle{definition}
\newtheorem{definition}{Definition}[section]
\newtheorem{example}{Example}[section]

\theoremstyle{remark}

\newtheorem{remark}{Remark}[section]

\numberwithin{equation}{section}

\newcommand{\ep}{\varepsilon}

 \DeclareMathOperator{\interior}{int}

   \DeclareMathOperator{\diam}{diam}

\begin{document}
\title[The Planck boundary]{The Planck boundary within the hyperspace of the circle of pseudo-arcs}
\author{Andr\'{a}s Domokos and Janusz R. Prajs}
\address{Department of Mathematics and Statistics,
California State University Sacramento, 6000 J Street, Sacramento, CA, 95819, USA}

\date{\today}



\begin{abstract}
In this paper we point out an interesting geometric structure of nonnegative metric curvature emerging from the hyperspaces of decomposable, non-locally connected homogeneous continua, where ``smooth" and ``non-smooth" partitions live together,  similarly to the macroscopic and the quantum realms of the Universe.  

\end{abstract}

\maketitle


\section{Introduction}
The main goal of this paper is to reveal interesting geometric  structures emerging from the study of hyperspaces of decomposable, non-locally connected homogeneous continua, in particular the hyperspaces of the circle of pseudo-arcs and the Menger curve of pseudo-arcs. 

These structures exhibit properties which, potentially,  could have applications to several disciplines. Therefore, we opted to take some highly abstract ideas from the area of topological continuum theory, and present them to scientists and mathematicians who are not necessarily involved in research  in continuum theory. 

The first part of the paper is an expository presentation of some foundational ideas of this theory. We blend together some definitions, properties and examples. In Section 3 we explain the duality between arcs and pseudo-arcs, which triggered many questions we try to answer in this paper. In Section 4 we review the definitions and properties of ample and filament subcontinua, which are key components of our new constructions. In Section 5 we explain the role of the Whitney map in generating the 
isometries on order-arcs, which are essential in establishing the nonnegative metric curvature of our hyperspaces. In Section 6 we prove a characterization of the circle of pseudo-arcs and show a method which generates a large family of circles of pseudo-arcs. In Section 7 we define and characterize the Planck boundary of the hyperspace of the circle of pseudo-arcs and prove the main theorem on its intrinsic metric geometry determined by the order-arcs. In Section 8 we propose some philosophical conclusions and questions related to our ideas. 
\smallskip

 All topological spaces in this paper are metric. A {\it continuum}
 is a non-empty, compact,  connected metric space. We will call  a continuum degenerate if it contains only one point. Compared to only visual observations, studying the structures of continua offers a deeper view on the variety of constructions ``one-piece" objects can have.
 \smallskip

In a compact space $X$ we introduce the {\it Hausdorff distance} or {\it Hausdorff metric} between non-empty closed subsets $K,L\subset X$ as
$$d_H (K,L) = \min \{ \ep > 0 \; : \; K \subset N (L, \ep) \; \text{and} \; L \subset N (K,\ep) \} \, ,$$
where
$$N(K, \ep) = \{ x \in X \; : \; d(x,k) \leq \ep \; \text{for some} \; k \in K \} \, $$
is the closed $\ep$-neighborhood of K.

 Following traditional continuum theory language, we refer to  continuous functions between topological spaces as {\it maps}. It is known that all  surjective maps $f: X \to Y$ between compact spaces are {\it quotient maps}, that is, the topology of $Y$ is uniquely  determined by the partition of $X$  into the point-inverses $f^{-1}(y)$. If the assignment $y\mapsto f^{-1}(y)$ is continuous with respect to the Hausdorff distance in $X$, we say that the partition ${\mathcal D} = \{ f^{-1} (y) : y \in Y \} $ is a {\it continuous decomposition} of $X$. The continuity of this decomposition is equivalent to the map $f:X\to Y$ being {\it open}, which means $f(U)$ is open in $Y$ for every open subset $U$ of $X$.

We call a continuum $X$ {\it homogeneous} if the group of self-homeomor-phisms of $X$ acts transitively on $X$, which means  for all $x, y \in X$ there exists a homeomorphism $h : X \to X$ such that $h(x) = y$. For example, the closed interval $[0,1] $ is not homogeneous, while the unit circle $S^1 \subset {\mathbb R}^2$ is. 
   
If $K$ is a non-empty, compact, connected subset of $X$, we say $K$ is a   {\it subcontinuum} of $X$. A commonly known subcontinuum, if exists, is the range of a one-to-one map $f: [0,1] \to X$, which we call an {\it arc}. 
 
For a continuum $X$ we define the {\it hyperspace} $C(X)$ as the collection of all subcontinua of $X$. The hyperspace $C(X)$ is a compact, connected  space with respect to the Hausdorff metric. In fact, $C(X)$ is also arcwise-connected, thanks to the existence of   {\it order-arcs} of subcontinua, which are defined as the ranges of maps $f: [0,1] \to C(X)$ such that $f(s) \subsetneqq f(t)$ whenever $0 \leq s < t \leq 1$.

As a simplistic, first view on a hyperspace of a homogeneous continuum $X$, we can imagine a  structure resembling a disk which has the degenerate subcontinua of $X$ (the set of  points of $X$)  forming the outer perimeter and the continuum $X$ being the center. Each point from the perimeter is connected by order-arcs to the center. As we move inward, the subcontinua (symbolized by points) gradually become larger and contribute to the change of the structure of the hyperspace.\\

\hspace*{3.2cm} \includegraphics[scale=0.4]{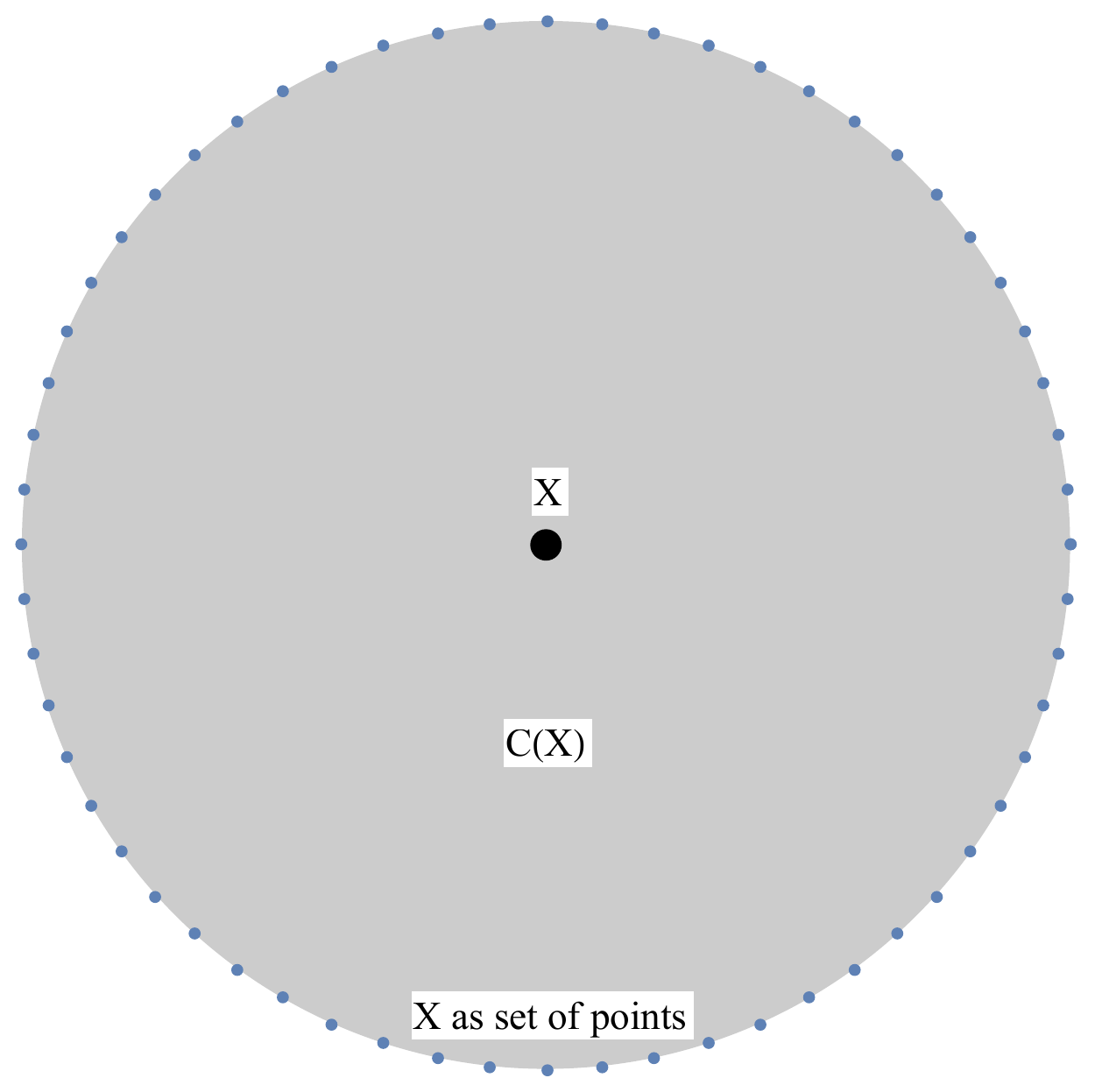}\\

To further gain some insight, we continue with two examples of 
hyperspaces. For similar examples we refer the reader to \cite{Nadler}. 

\begin{example} In this first example we consider the non-homogeneous continuum $X = [0,1]$. Consider the triangular region  $T$ with vertices at $(0,0)$, $(1,0)$ and $(0.5,1)$. Each subcontinuum is a point or a closed interval $[a,b] \subset [0,1]$. We can define a homeomorphism $H : C(X) \to T$, as
$$H([a,b]) = \Bigl( \, \frac{a+b}{2} \, , \, b-a \, \Bigr) \, .$$	
Of course, if $a=b$ the above formula gives the value for a degenerate subcontinuum: $H(\{a\}) = (a,0)$.\\
\hspace*{3.7cm} \includegraphics[scale=0.4]{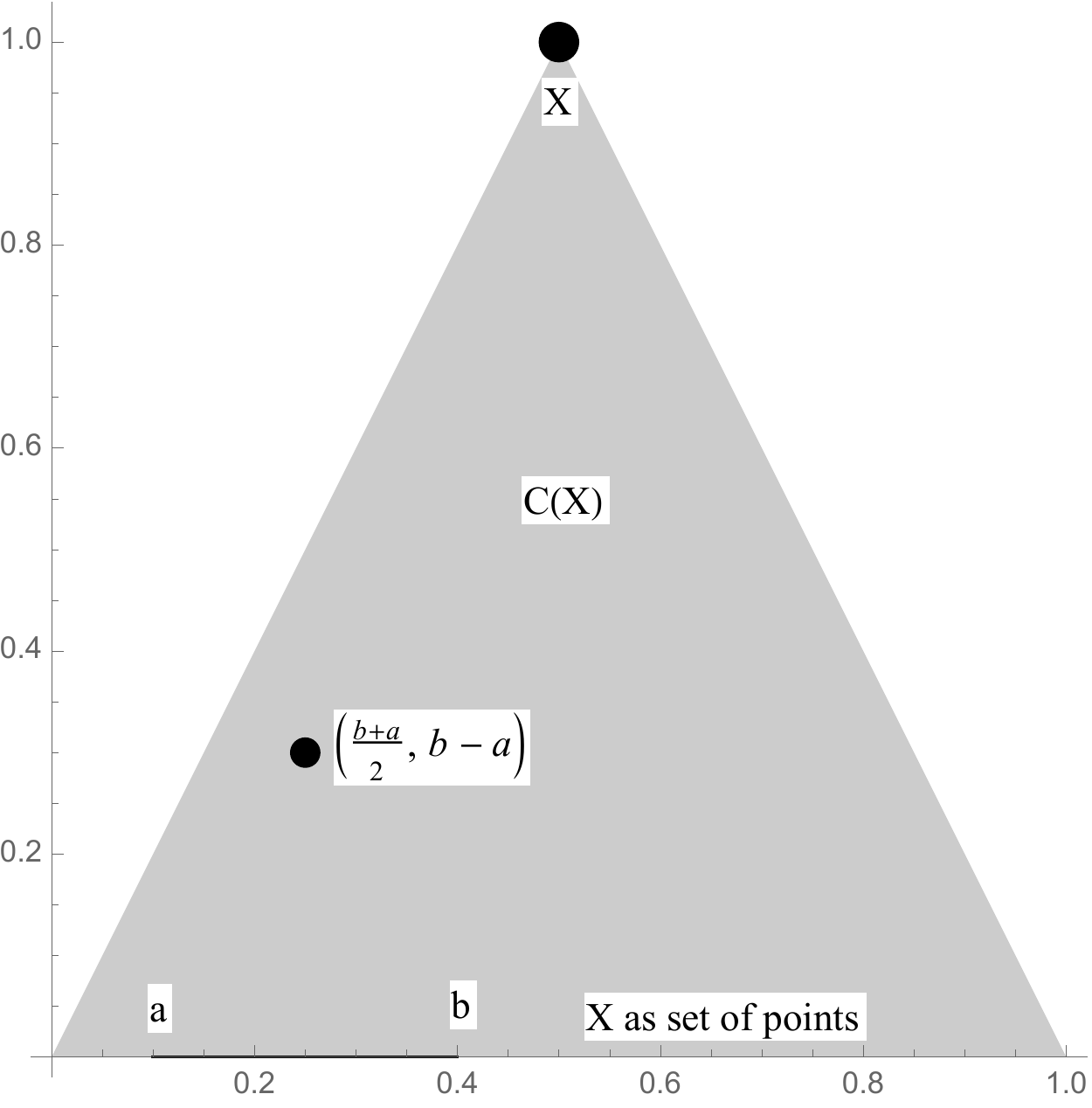} \\
It is easy to check that $H$ is one-to-one, onto and continuous. The compactness of the C(X) implies that the inverse function is continuous, too, therefore $H$ is a homeomorphism.
\end{example}
 
\begin{example}
As the second example, consider the unit circle $S^1$ in the plane, which is a homogeneous continuum. Each subcontinuum is a point or an arc of the circle $[\alpha, \beta]$, where $\alpha$ and $\beta$ are the angles defining the starting and ending points of the arc.
We can define a homeomorphism $H : C(S^1) \to D$, where $D$ is the closed unit disk, as
$$H([\alpha,\beta]) = \Bigl(1-\frac{\beta-\alpha}{2\pi} \Bigr) \; e^{\frac{\beta+\alpha}{2}i} \, .$$

We can visualize this homeomorphism with the following picture.\\

\hspace*{3.2cm} \includegraphics[scale=0.4]{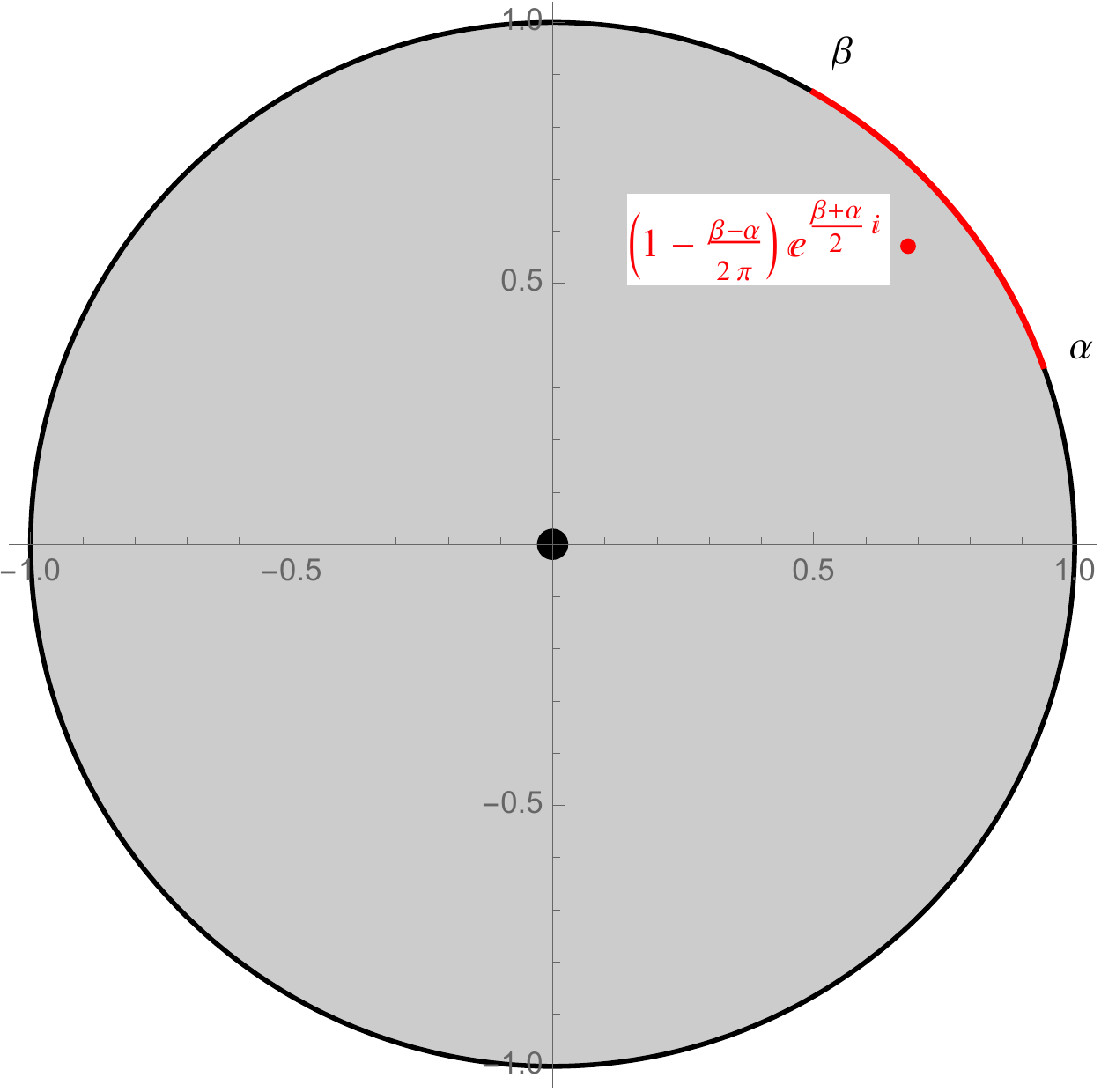} \\	

Also, we can show a graph with some order-arcs emanating from the point $(1,0)$.\\

\hspace*{3.2cm} \includegraphics[scale=0.4]{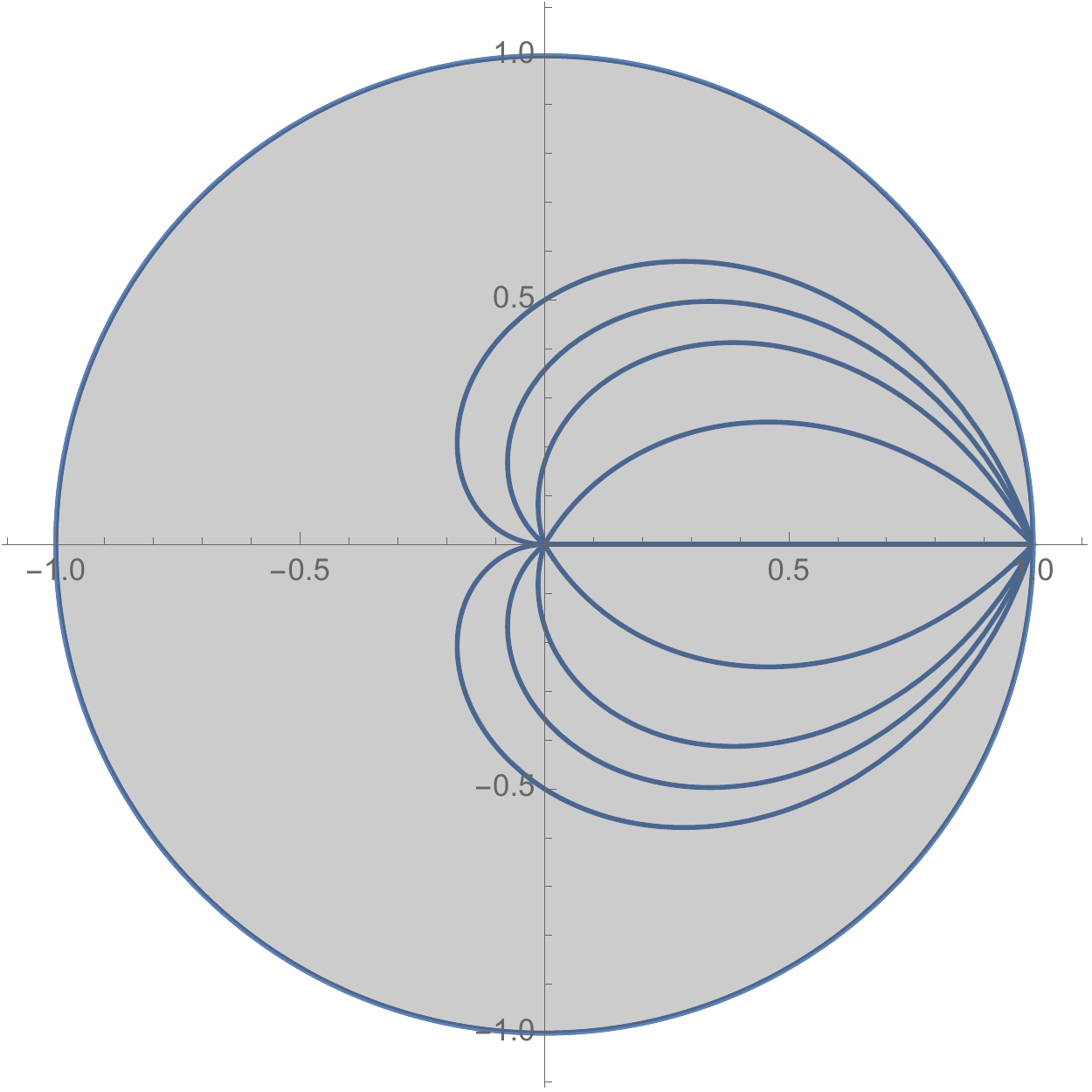}\\

The two extremal order-arcs bound the local accessibility region at the point (1,0). Similar situations appear at every point of the disk. This suggests a certain type of length metric geometry for $C(S^1)$, defined by connectedness through concatenations of order-arcs. As we will see later, in general situations, the order-arcs connect the ``smooth" and ``non-smooth" sides of the hyperspace. 

A similar construction can be done by embedding $C(S^1)$ in ${\mathbb R}^3$ as the cone over the unit circle.
\end{example}


\section{Local connectedness, decomposability and indecomposability}

A continuum $X$ is {\it locally connected at} $x \in X$ if every neighborhood of $x$ includes a connected neighborhood of $x$. The continuum $X$ is said to be {\it locally connected} if it is locally connected at every point. If  $X$ is not locally connected at any point, we say $X$ is {\it nowhere locally connected}. Each homogeneous continuum  is  either locally connected, or nowhere locally connected. The unit circle is locally connected,  while the pseudo-arc and circle of pseudo-arcs, defined in the following sections, are nowhere locally connected.

A continuum $X$ is said to be {\it decomposable} if there exist $A, B \in C(X)$ such that 
$A \neq X\neq B$ and $X = A \cup B$.  For instance, an arc is the union of two smaller arcs, a circle is the union of two arcs, a $2$-sphere is the union two hemispheres etc. In general, every non-degenerate locally connected continuum is decomposable.  If every non-degenerate subcontinuum of $X$ is decomposable, we say $X$ is {\it hereditarily decomposable}.

\begin{example}\label{example-cone-Cantor}
For visual insight we include a graph of the cone over the ternary Cantor set, called the {\it Cantor fan}. \\
\hspace*{4cm} \includegraphics[scale=0.4]{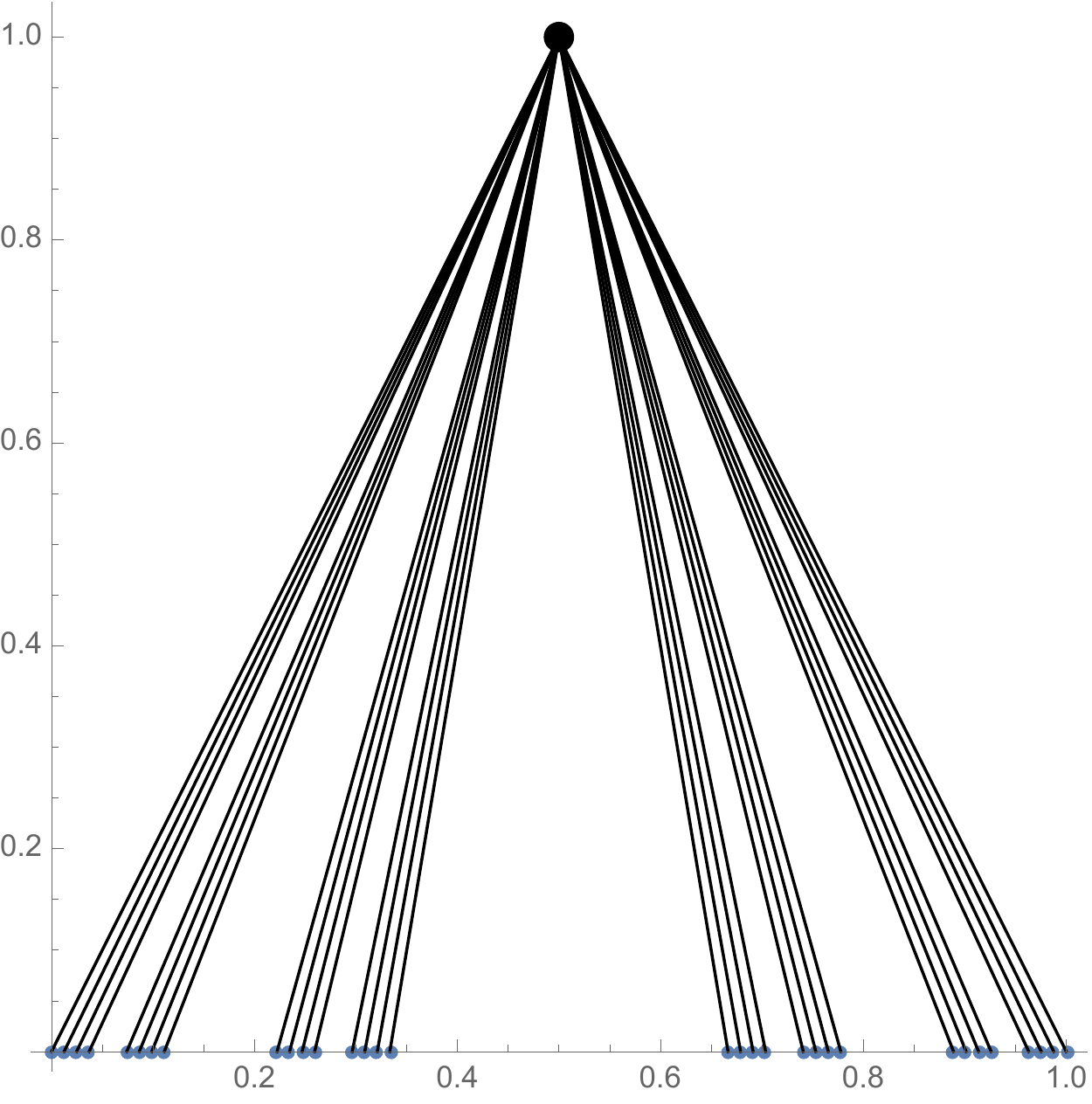} \\
\indent The Cantor fan is non-homogeneous, hereditarily decomposable, and it is  locally connected only at the vertex.	
\end{example}

A continuum is {\it indecomposable} if it is not decomposable. Indecomposable continua are harder to visualize and more counterintuitive than the decomposable ones. The study of indecomposable continua began around 1912 with an example now known as the {\it bucket handle continuum} attributed to various authors. Descriptions, sketches and properties of this continuum are easily avaiable in the literature and in the Internet. In the 1960s this continuum was rediscovered in the study of dynamics under the name {\it Smale's horseshoe}.
\smallskip


 A continuum is {\it hereditarily indecomposable} if each of its subcontinua  is indecomposable. The first example of a hereditarily indecomposable continuum was given by B. Knaster in 1922 \cite{Knaster-22}.
 \smallskip
 
 Many constructions of  decomposable, indecomposable and hereditarily indecomposable continua use the following theorem.

\begin{theorem}  \label{nested}
Consider a metric space $M$ and a sequence $\{ X_n \}$ of non-empty compact, connected subsets of $M$ such that $X_{n+1} \subset X_n$ for all $n \in {\mathbb N}$. Then
$$X = \bigcap_{n=1}^{\infty} X_n$$
is a non-empty, compact, connected subset of $M$. Moreover,
$$\lim_{n \to \infty} d_H (X_n , X) = 0 \, .$$
\end{theorem}


\section{The Duality between an Arc and the Pseudo-Arc}\label{Section-Duality}\label{Section 3}

The significance of the topological structure of an {\it arc}, a homeomorphic copy of the unit segment
$[0,1]$, cannot be overestimated. It is a fundamental topological building block for the real number line, linear spaces, manifolds, and is used in almost all areas of mathematics and its applications. The {\it pseudo-arc} is mostly unknown beyond the area of topology. In this section we compare the properties of these two special topological objects.

A continuum $X$ is called {\it chainable} if for every $\varepsilon > 0$ there is a finite open cover  $C_1, C_2, ... , C_n$ of $X$ such that for all $k,l\in \{1, .... , n\}$ $\diam C_k< \ep$, and $C_k \cap C_l \neq \emptyset$ if and only if $|k-l| \leq 1$. The collection ${\mathcal C}=\{C_1, C_2, ... , C_n\}$  is called an $\ep$-{\it chain}, and the sets $C_i$ are called {\it links}. We write $\cup \, {\mathcal C} = C_1 \cup ... \cup C_n$.
\smallskip

The foundational work on the existence and properties of the  pseudo-arc was done by Knaster \cite{Knaster-22},   Moise \cite{Moise} and     Bing \cite{Bing-48, Bing-51}. Here we offer an abrieviated version of the construction of the pseudo-arc.

\begin{example} (The pseudo-arc)

 Consider an $\ep$-chain ${\mathcal C}$ and a $\delta$-chain ${\mathcal D}$, both composed of open topological discs in the plane. (We think $\delta$ is substantially less than $\ep$.) 
  
 We say that the chain ${\mathcal D}$ is crooked in ${\mathcal C}$ if the following two conditions hold:\\
 (1) For each link $D_i$ of $\mathcal D$, there is a link $C_k$ of $\mathcal C$ such that $\overline{D_i} \subset C_k$.\\
 (2) For any pair of indices $k$ and $m$ such that $k \leq m-3$ and any $k \leq i < j \leq m$ such that $D_i \cap C_k \neq \emptyset$ and $D_j \cap C_m \neq \emptyset$, there exist indices $i_*$ and $i^*$ such that $i < i_* < i^* < j$ and 
 $\overline{D_{i_*}} \subset C_{m-1}$ and  $\overline{D_{i^*}} \subset C_{k+1}$.
 
 For a sample of the above definition of crookedness, see the following picture. Note that the links are the regions (not the curves) enclosed by their boundaries.\\
 
 \hspace*{2cm} \includegraphics[scale=0.3]{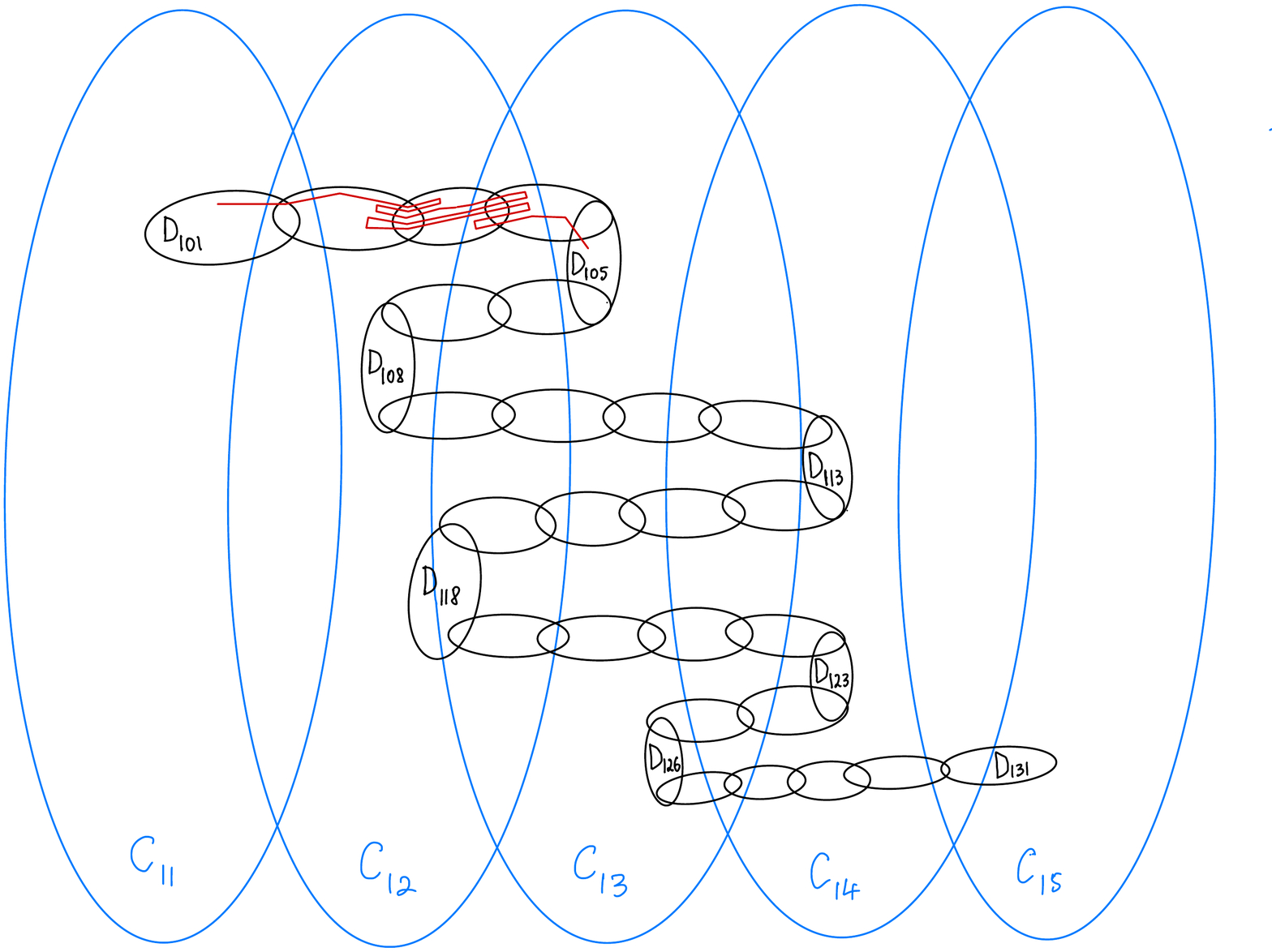} \\	

\begin{definition}
Let $x \neq y \in {\mathbb R}^2$. Consider a sequence of chains $\{ {\mathcal C}^n \}$ such that for each  $n \in {\mathbb N}$:\\
(a)   $x$ is in the  the first  and $y$ in the last link of $ {\mathcal C}^n $.\\
(b)  ${\mathcal C}^n$ is a $2^{-n}$-chain.\\
(c)  ${\mathcal C}^{n+1}$ is crooked in ${\mathcal C}^n$. 
\smallskip

Then, we define the pseudo-arc connecting $x$ and $y$ as
$${\mathcal P} = \bigcap_{n=1}^{\infty} \overline{\bigcup{\mathcal C}^n}\, .$$

The set $\mathcal{P}$ is a continuum in view of Theorem \ref{nested}. As it was proved in the papers mentioned above this example, the pseudo-arc is  hereditarily indecomposable. In particular, it contains no arcs. 
	
\end{definition}

\end{example}

{\bf Similarities between  arcs and pseudo-arcs.}
\smallskip

Both are {\it chainable}.
\smallskip

Both are {\it topologically unique}, which means that any two arcs are homeomorphic and any two pseudo-arcs are homeomorphic. An arc is unique by definition. A major breakthrough was the result by Bing \cite{Bing-51} showing that each two hereditarily indecomposable, chainable continua are homeomorphic. This characterization of the pseudo-arc connected Knaster's example \cite{Knaster-22} with Moise's example \cite{Moise} as topologically the same.
\smallskip

Both are {\it hereditarily equivalent}, that is, every non-degenerate subcontinuum is homeomorphic to the whole continuum. This is why Moise selected the name {\it pseudo-arc} for his example. Whether these two are the only hereditarily equivalent continua is an intriguing question in continuum theory, which has been open for many decades.
\smallskip

Both are {\it irreducible}, that is, contain a pair of points such that no proper subcontinuum contains these two points.
\smallskip

{\bf Contrasting differences between arcs and pseudo-arcs.}
\smallskip

An arc is locally connected while the pseudo-arc is nowhere locally connected.
\smallskip

An arc is hereditarily decomposable while the pseudo-arc is hereditarily indecomposable.
\smallskip

An arc is irreducible between exactly one pair of points while the pseudo-arc is irreducible between uncountably many pairs of points.
\smallskip

The pseudo-arc is homogeneous \cite{Bing-48} while an arc is not.
\smallskip

Since the pseudo-arc is homogeneous, it has rich spaces of self-maps and self-homeo\-morphisms.   Yet, the pseudo-arc is {\it pseudo-homotopically rigid} at any self-homeomorphism \cite{Illanes}. ({\it Pseudo-homotopy } is similar to homotopy, but the parameter space can be any continuum and not necessarily an arc.) On the other hand, an arc is contractible, or equivalently, the space of self-maps of an arc is arcwise connected.
\smallskip

Most continua, in some sense, are pseudo-arcs. Let $X$ be an $n$-manifold, $n>1$,  or a Hilbert cube manifold. Then the hyperspace of subcontinua, $C(X)$, contains uncountably many mutually non-homeo\-morphic continua. Yet, the copies of a single continuum, the pseudo-arc, form a dense $G_{\delta}$ set in $C(X)$ \cite{Bing-51}. That means the copies of all continua other than the pseudo-arcs, together, including arcs, are a ``thin" first Baire category collection.
\smallskip

In any metric space $(M,d)$ the convergence of a sequence  $P_n$  of pseudo-arcs to a pseudo-arc $P_0$ with respect to the Hausdorff metric implies the {\it homeomorphic convergence}, which means there are homeomorphisms $h_n:P_n \to P_0$ such that
$\lim_n(\sup\{ d(x,h_n(x))\,:\,x\in P_n\}) = 0$. As it can be seen on the picture below, an arc does not have this property. In fact, among all continua, only the pseudo-arc possesses  this extraordinary property \cite{Lewis-98}.\\

 \hspace*{3cm} \includegraphics[scale=0.3]{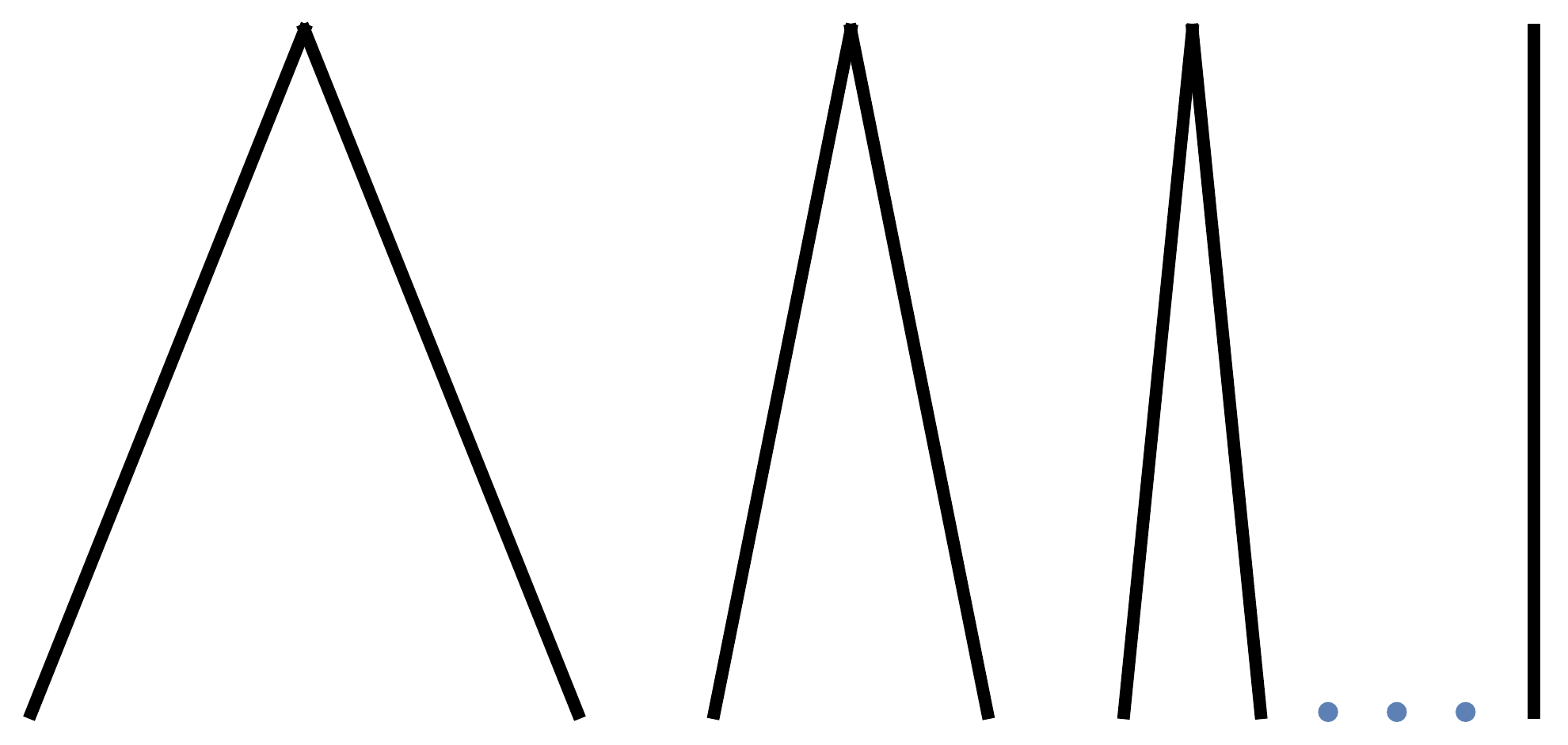} \\

 There are numerous surprising constructions  involving pseudo-arcs \cite{B-J-59} \cite{LewisWalsh} \cite{Lewis-85} \cite{Prajs1} \cite{Prajs2} \cite{Prajs3} \cite{Prajs4}, which are very different from constructions involving arcs. Many of these constructions are continuous decompositions into pseudo-arcs. Any such decomposition must be {\it completely regular}, that is, it is not only continuous with respect to the Hausdorff distance, but it must also be continuous with respect to the homeomorphic convergence.
\smallskip
\smallskip
\smallskip

In view of these properties, the pseudo-arc has a very special position  in metric topology and in the fabric of important spaces such as manifolds. These properties have been a repeating topic of conversations between the authors of this paper. A fundamental question arises from these conversations.\\

{\bf Question.}  Could the pseudo-arc and constructions involving the pseudo-arc be used as models outside topology? Outside mathematics? \\

The main motivation for this paper was triggered by this question.



\section{Filament and ample subcontinua}

In this section we follow the definitions and properties of filament and ample continua as developed in \cite{Pr-W-06}.


A subcontinuum $K \in C(X)$ is called  {\it filament} if there exists an open set $U$ such that $K \subset U$ and the component of $U$ containing $K$ has empty interior.
As we can see in Example \ref{example-cone-Cantor}, in the Cantor fan every subcontinuum, which does not contain the vertex, is filament. In the pseudo-arc, as in every indecomposable continuum, every proper subcontinuum is filament.

A subcontinuum $K \in C(X)$ is called {\it ample} if for every open set $U$ with $K \subset U$, there exists $L \in C(X)$ such that $K \subset \interior(L) \subset L \subset U$. Every subcontinuum of a locally connected continuum is ample. For instance, every subcontinuum of a compact connected manifold is ample.  In the Cantor fan every subcontinuum containing the vertex is ample.

In a homogeneous continuum $X$, a subcontinuum is ample if and only if it is not filament. Let us denote by $\mathcal{F}$ the set of filament subcontinua and call it the filament portion of $C(X)$. Similarly, we denote by $\mathcal{A}$ the set of ample subcontinua and call it the ample portion of  $C(X)$. Consequently, for a homogeneous continuum $X$ we have $C(X) = \mathcal{F} \cup \mathcal{A}$ and $\mathcal{F} \cap \mathcal{A} = \emptyset$. 

The size of $\mathcal{F}$ shows where $X$ stands between  local connectedness  on one end and indecomposability on the other end.
A homogeneous continuum $X$ is locally connected if and only if $\mathcal{F} = \emptyset$, and it is indecomposable if and only if
$\mathcal{F} = C(X) \setminus \{X\}$.

Regarding the properties of filament and ample subcontinua of a homogeneous continuum, we mention the following.

If $K \in C(X)$, $L \in {\mathcal F}$ and $K \subset L$, then $K \in {\mathcal F}$.

 
For each $K \in {\mathcal F}$ there exists $\delta > 0$ such that $L \in {\mathcal F}$ for every $L \in C(X)$ with $d_H (K,L) < \delta$. Thus ${\mathcal F}$ is an open subset of $C(X)$. Also, if ${\mathcal F}$ is nonempty, it contains all degenerate subcontinua of $X$, which form a homeomorphic copy of $X$ embedded in ${\mathcal F}$. Order-arcs in ${\mathcal F}$ connect every element of ${\mathcal F}$ to a degenerate subcontinuum. Therefore, 
${\mathcal F}$ is a connected and open subset of $C(X)$. 

If $K \in {\mathcal A}$, $L \in C(X)$ and $K \subset L$, then $L \in {\mathcal A}$. This shows that any $K \in {\mathcal A}$ can be connected  to $X$ by an order-arc in ${\mathcal A}$.

It is known that $K \in {\mathcal A}$ if and only if $C(X)$ is locally connected at $K$. In addition, by Proposition 2.16 \cite{Pr-W-06}, $\mathcal A$ is an absolute retract and therefore it is locally contractible. 

In the next theorem we summarize the properties of the filament and ample portions in the context relevant to the hyperspace of the circle of pseudo-arcs discussed in Section \ref{Section-circle of pseudo-arcs}.  

\begin{theorem}\label{prop-ample-filament} In a decomposable and nowhere locally connected homogeneous continuum $X$ the following properties hold:\\
(a) The ample portion ${\mathcal A}$ is a proper subcontinuum of $C(X)$, and it is locally arcwise-connected, contractible and locally contractible.\\
(b) Every $A \in {\mathcal A}$ contains a minimal ample subcontinuum, which  either is indecomposable or it has empty interior.\\
(c) The filament portion ${\mathcal F}$ is non-empty, open, connected and nowhere locally connected. \\
(d) $C(X) = \mathcal{F} \cup \mathcal{A}$ and $\mathcal{F} \cap \mathcal{A} = \emptyset$.  
\end{theorem}



\section{Whitney maps and distances}\label{Section-Whitney}

In this section we use the usual topology notation of $2^X$ for the collection of all non-empty,  closed (hence compact) subsets of a continuum $X$. The collection $2^X$ endowed with the Hausdorff distance is also a continuum. We need the space $2^X$ because the definition of our Whitney map requires that non-empty intersections stay in the space. Once we obtain the Whitney distance on $2^X$, we will restrict it to $C(X)$.

\begin{definition}\label{def-Whitney}
A map $\mu : 2^X \to [0,+\infty)$	 is called a Whitney map if it satisfies the following conditions:\\
(a) $\mu (\{x\}) = 0$, for all $x \in X$.\\
(b) If $A,B \in 2^X$ and $A \subsetneqq B$, then $\mu (A) < \mu (B)$.\\
(c) If $A, B \in 2^X$ and $A \cap B \neq \emptyset$, then
$$\mu(A\cup B) \leq \mu(A) + \mu(B) - \mu(A\cap B) .$$
\end{definition}

\begin{remark} 
The usual definition of a Whitney map includes only the conditions (a) and (b). However, the extra condition (c) allows the definition of a metric which, when restricted to $C(X)$, provides an interesting metric geometry structure.
\end{remark} 

Definition \ref{def-Whitney} is equivalent to the definition of a Whitney map from \cite{Charat}. The difference is that condition (c) is replaced by the following:

(c') If $A',B' \in 2^X$ and $A' \subset B'$, then
$$\mu(B'\cup C') - \mu(A'\cup C') \leq \mu(B') - \mu(A') \, , \; \text{for all} \; C' \in 2^X \, .$$

We opted for this formulation of (c), because it relates a finitely sub-additive measure on some collections of compact subsets of $X$ to a metric and its geometry on $C(X)$.

\begin{lemma}
The conditions (a), (b) and (c) from Definition \ref{def-Whitney}  are equivalent to (a), (b) and (c').	
\end{lemma}

\begin{proof} Assume that (a), (b) and (c) hold. Consider $A',B' , C' \in 2^X$ such that $A' \subset B'$.  We use (c) for $B= B'$ and $A= A' \cup C'$. Note that $(A'\cup C')\cap B'=(A'\cap B')\cup(C'\cap B')=A' \cup (B'\cap C')$. Thus,
\begin{align*}
\mu (B' \cup C')&=\mu(A'\cup C' \cup B') \\
& \leq \mu (A' \cup C') + \mu (B') - \mu (A' \cup (B' \cap C'))
\\ &\leq \mu (A' \cup C') + \mu (B') - \mu (A') \, ,
\end{align*}
from which (c') follows.
\smallskip

Now assume that (a), (b) and (c') hold. Consider $A,B \in 2^X$ such that $A\cap B\neq\emptyset$. If $A \subset B$, then (c) holds. Assume that $A \not\subset B$. We use (c') for $A' = A \cap B$,\, $B' = B$
and 
%
%
$C' = \overline{A\setminus B}$. Note that $A\cup B = B \cup C'$ and $A=(A\cap B)\cup C'$.
Thus,
\begin{align*}
\mu(A\cup B) - \mu(A)&=\mu(B\cup C')- \mu((A\cap B)\cup C') \\
&\leq \mu(B) - \mu(A\cap B) ,
\end{align*}
from which (c) follows.
\end{proof}

\begin{prop}\label{prop-Whitney-exists} \cite{Charat}
For every continuum $X$, a Whitney map in the sense of Definition \ref{def-Whitney} exists.	
\end{prop}

\begin{proof}
As mentioned in \cite{Charat}, the construction 	(0.50.2) in \cite{Nadler} provides a Whitney map with the properties (a), (b) and (c'). For the completeness of the presentation, we show the construction of this Whitney map. 

Consider a countable dense subset $\{x_1, x_2, ... \}$ of $X$. For each $n \in \mathbb N$ define the functions
$$f_n : X \to [0,1] \, , \; f_n (x) = \frac{1}{1+d(x_n, x)} \, ,$$
and 
$$\mu_n : 2^X \to [0, +\infty) \, , \; \mu_n( A ) = \diam (f_n (A)) \, .$$
Now define the Whitney map as
$$\mu:2^X \to [0,+\infty) \, , \; \mu(A) = \sum_{n=1}^{\infty} \frac{\mu_n (A)}{2^n} \, .$$
\end{proof}

\begin{definition}\label{def-Whitney-metric} \cite{Charat} For all $A,B \in 2^X$ define the Whitney distance as 
$$d_{\mu} (A,B) = \max \{ \mu(A\cup B) - \mu (A), \mu (A\cup B) - \mu(B) \} \, .$$
\end{definition}

The next proposition lists the properties of $d_{\mu}$.

\begin{prop}\label{prop-d_mu} \cite{Charat} The following statements hold:\\
(a) $d_{\mu}$ is a metric on $2^X$.\\
(b) $d_{\mu}$ is equivalent to the Hausdorff metric $d_H$.\\
(c) $d_{\mu}(A , \{x\}) = \mu(A)$ for all $x \in A \in 2^X$.\\
(d) For any order arc ${\mathcal O}$ in $2^X$, $\left. \mu \right|_{\mathcal O} : {\mathcal O} \to [0, \mu (X)]$ is an isometry. 
\end{prop}

\begin{remark}
For the proofs of (c) and (d), observe that  from Definition \ref{def-Whitney-metric} follows that if $A \subset B$, then
$$d_{\mu} (A,B) = \mu (B) - \mu (A) \, .$$	
\end{remark}

\begin{remark}
The equivalence of $d_{\mu}$ and $d_H$ implies that 	$C(X)$ is a continuum with respect to $d_{\mu}$.
\end{remark}

As it was stated in Section \ref{Section 3}, many spaces and continua admit continuous decompositions into pseudo-arcs. Among them are the plane, all compact connected surfaces with and without boundary, the Sierpi\'nski `carpet' curve, and the Menger curve. We end this section with   showing how to modify these decompositions so that all  the members of the decomposition are  on the same Whitney level.

\begin{prop}\label{propDecomp}
 Let $X$ be a continuum with a continuous decomposition $\mathcal{D}$ into pseudo-arcs and $\mu: C(X) \to [0,\mu(X)]$ be a Whitney map. Then there exists a continuous decomposition $\mathcal{D}_0$ of $X$ into pseudo-arcs, which is finer than $\mathcal{D}$, such that $\mu(D)$ is constant for all $D\in \mathcal{D}_0$.
\end{prop}

\begin{proof}
First, notice that for a pseudo-arc $P$ in $X$ with $\mu(P)= s$, and for all $t\in[0,s]$ and  $x\in P$ there is a unique continuum $P(x,t)\subset P$ such that $x\in P(x,t)$ and $\mu(P(x,t))=t$. This is a consequence of the heredirary indecomposability of $P$. 

Since $X$ is compact and $\mathcal{D}$ is continuous, $\mu(D)$  cannot be arbitrarily small for $D\in \mathcal{D}$. Let $r =\min\{\mu(D)\,:\, D\in \mathcal{D}\}$ and $t_0\in (0,r]$. Define 
$$\mathcal{D}_0=\{ D_0 \in C(X)\,:\, \mu(D_0) =t_0\,\,\mathrm { and } \,\, D_0 \subset  D \,\, \mathrm{ for \,\, some } \,\, D\in \mathcal{D}\}.$$
 By the hereditary indecomposability of pseudo-arcs and properties of Whitney maps,  $\mathcal{D}_0$ is a continuous partition of $X$ into pseudo-arcs.
\end{proof}


\section{The circle of pseudo-arcs}

A subcontinuum $K \subset X$ is {\it  terminal} in $X$, if for any continuum $L \subset X$ with $K \cap L \neq \emptyset$ we have either $K \subset L$ or $L \subset K$. A map $f:X\to Y$ between continua $X$ and $Y$ is called {\it atomic} if  all point-inverses $f^{-1}(y)$ for $y\in Y$ are terminal subcontinua of $X$.  

\begin{example} (The circle of pseudo-arcs)

In 1955, Jones \cite{Jones-55} showed that every homogeneous, decomposable plane continuum has a continuous decomposition into mutually homeomorphic, homogeneous, non-separating plane continua, such that the quotient space  defined by the decomposition is a simple closed curve. Triggered by this result, Bing and Jones \cite{B-J-59} constructed the {\it circle of pseudo-arcs}, denoted  here by  $\Psi$, which is a homogenous, decomposable and nowhere locally-connected continuum in the plane 
admitting a continuous decomposition into terminal pseudo-arcs. Therefore, we can assume having an open atomic map $\Pi : \Psi \to S^1$ such that for all $\alpha \in [0, 2\pi)$ there exists a terminal pseudo-arc $P_{\alpha}$ such that
$$\Pi^{-1} (e^{i \alpha}) = P_{\alpha} \, ,$$
and
\begin{equation}\label{eq:Jones-decomposition}
\Psi = \bigcup_{\alpha \in [0,2\pi)} P_{\alpha} \, .
\end{equation}
\end{example}

We call a continuum a circle of pseudo-arcs if it is homeomorphic to the original example $\Psi$ of Bing and Jones \cite{B-J-59}.

A {\it triod}\footnote{A {\it triod}  should not be confused with the {\it simple triod}, which is the union of three distinct arcs emanating from a common end-point.} is the union $T$ of  three continua $A$, $B$ and $C$ such that $K=A\cap B=B\cap C =A\cap C$ is a proper subcontinuum of each of $A$, $B$ and $C$.  

\hspace*{3cm}\includegraphics[scale=0.35]{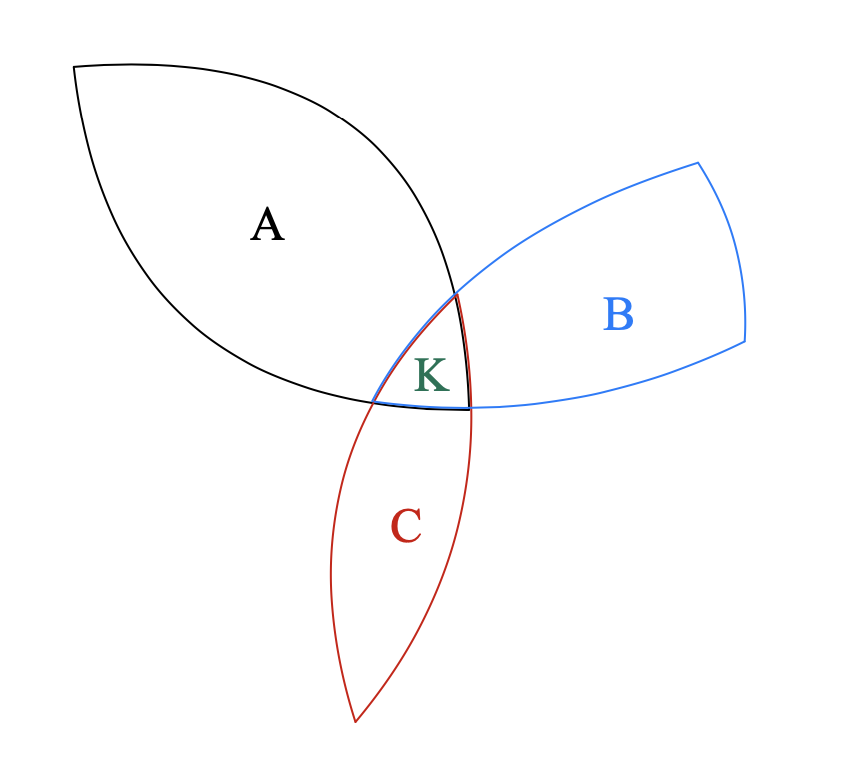} 

A continuum $X$ is said to be {\it atriodic} if it contains no triods. An arc, a circle, the pseudo-arc and the circle of pseudo-arcs are examples of atriodic continua.

The next theorem gives a characterization of the circles of pseudo-arcs. Compared to the characterization  in \cite{B-J-59}, we do not use {\it circular chainability}. We include the proof because we did not find in the literature a proper reference to this theorem.


\begin{theorem}\label{uniqueness} 
A continuum $X$ is a circle of pseudo-arcs if and only if $X$ admits an open atomic map $f:X\to S$ onto a simple closed curve $S$ such that all point-inverses, $f^{-1}(y)$, are pseudo-arcs.

In particular, the circle of pseudo-arcs is topologically unique.
\end{theorem}

\begin{proof}
It is known that the circle of pseudo-arcs from \cite{B-J-59} admits an open, atomic map onto a circle with pseudo-arcs as fibers. 

Let $X$ be a continuum with an open atomic map $f:X\to S^1$ onto the unit circle in the plane. 
Using the atriodicity of $S^1$ and of the pseudo-arc, and the atomicity of $f$ one can easily show that $X$ is atriodic. 

Define $A=\{(x,y)\in S^1\,:\,y\geq 0\}$ and $B=\{(x,y)\in S^1\,:\,y\leq 0\}$. Then $f^{-1}(A)$ and $f^{-1}(B)$ are atriodic folders of chainable continua in the sense of \cite{HMP}. Therefore, they are chainable by Proposition 6 of \cite{HMP}. Using Theorem 10 of  \cite{B-J-59}, $X=f^{-1}(A) \cup f^{-1}(B)$ can be mapped homeomorphically onto the circle of pseudo-arcs of  \cite{B-J-59}.
\end{proof}

In \cite{Prajs1} it was shown that each compact, connected $2$-manifold $M$ admits  a continuous decomposition into pseudo-arcs such that the quotient space is homeomorphic to $M$. The next theorem offers a large collection of circles of pseudo-arcs that can be derived from these decompositions.

\begin{theorem}\label{pseudo-circlesinsphere} 
Let $\mathcal{D}$ be a continuous decomposition of the unit $2$-sphere $S^2$ into pseudo-arcs with the quotient map $q:S^2 \to S^2$  of $\mathcal{D}$. Then, for every simple closed curve $Z\subset S^2$, the set $q^{-1}(Z)$ is a circle of pseudo-arcs.
\end{theorem}

\begin{proof}
Let $q:S^2 \to  S^2$ be the quotient map from the assumption and $Z\subset S^2$ be a simple closed curve. Then $A=q^{-1}(Z)$ admits a continuous decomposition into pseudo-arcs $q^{-1}(y)$ for $y\in Z$.
    
By Theorem \ref{uniqueness} it suffices to show that each of the pseudo-arcs $q^{-1}(y)\subset A$  is  terminal in $A$. First, since $q$ is open and $q(A)=Z$ has no interior, $A$ has no interior either. Since $Z$ is the common boundary of  exactly two complementary domains $R_1$ and $R_2$ in $S^2$, and the decomposition of $S^2$ into the pseudo-arcs $q^{-1}(y)$ is continuous, it follows that $A=q^{-1}(Z)$ is the common boundary in  $S^2$ of exactly two complementary domains
$\widehat{R}_1=q^{-1}(R_1)$ and  $\widehat{R}_2= q^{-1}(R_2)$ in $S^2$. 

Suppose that there exists $\alpha \in Z$ such that $q^{-1}(\alpha)$ is not terminal in $A$. Let $K$ be a continuum in $A$ such that $q^{-1}(\alpha)\cap K\not= \emptyset\not= K\setminus q^{-1}(\alpha)$. Let $L$ be a component of $K\cap q^{-1}(\alpha)$. Given a Whitney map $\mu :2^A \to [0,1]$, we have
 $ \mu(L)<\mu(f^{-1}(\alpha))$. By the continuity  of the decomposition  of $A$ into the sets $q^{-1}(y)$, $y\in S$, and by the continuity of $\mu$, we can slightly enlarge  $L$ to a continuum $M \supset L$ such that the Whitney size $\mu(M)$ prevents $M$ from containing any fiber $q^{-1}(y)$ it intersects,  but $M\cap q^{-1}(\beta)\not=\emptyset $ for some $\beta \in Z$ different from $\alpha$. Since $Z$ is a simple closed curve, there are arcs $P$ and $Q$ in $Z$, both having $\alpha $ and $\beta$ as end-points, such that 
 $P\cup Q = Z$ and $q(M) \subset P$. By the known properties of separation of the sphere $S^2$, the continuum $B = M \cup q^{-1}(Q)$ separates $S^2$ between  any point $r_1\in \widehat{R}_1$ and any point  $r_2\in \widehat{R}_2$. 
  
Let $W_1$ and $W_2$ be the components of $S^2\setminus B$ containing $\widehat{R}_1$ and $\widehat{R}_2$, respectively. Fix a $\gamma \in P$ with $\alpha \not= \gamma \not= \beta$. Then neither $M$ nor  $q^{-1}(Q)$  contains $q^{-1}(\gamma)$. Therefore, there exists $p\in q^{-1}(\gamma) \setminus B$. Since $W_1\cap W_2=\emptyset$, the point $p$ cannot belong to both $W_1$ and $W_2$. Suppose $p \notin W_1$. Note that there is  a sequence $\{x_n\}$ in $R_1$ converging to 
$\gamma$ because $S$ is the common boundary of $R_1$ and $R_2$.  However, the sequence of pseudo-arcs $q^{-1}(x_n)$ cannot have the point $p\in q^{-1}(\gamma)$ in its limit set because because all  $q^{-1}(x_n)$ are separated from $p$ by $B$. Hence, the decomposition of $S^2$ into the sets $q^{-1}(y)$ is not continuous, a contradiction.
\end{proof}

\begin{remark}
One can show that for the map $q:S^2 \to S^2$, as in  Theorem \ref{pseudo-circlesinsphere}, and for every topologically one-dimensional continuum $X$ in $S^2$,
the map  $q|_{q^{-1}(X)} : q^{-1}(X) \to X$ is atomic. The argument is significantly more complicated, and this last claim goes beyond the scope of this paper.
\end{remark}


In 2016, Hoehn and Oversteegen \cite{H-O-16} showed that the only homogeneous plane continua (up to homeomorphism equivalence) are the unit circle, the pseudo-arc and the circle of pseudo-arcs.

For other non-planar homogeneous continua, constructed, based on the Jones decomposition theorem, similarly to the circle of pseudo-arcs, we can refer to the papers of Hagopian and Rogers \cite{H-R-77} and Lewis \cite{Lewis-83, Lewis-85}.


\section{The hyperspace of the circle of pseudo-arcs}\label{Section-circle of pseudo-arcs}

In this section we consider a circle of pseudo-arcs $\Psi$ which admits the decomposition \eqref{eq:Jones-decomposition}. As mentioned earlier, the circle of pseudo-arcs is a homogenous, decomposable and nowhere locally-connected continuum in the plane.

\begin{definition} We define the Planck boundary of the hyperspace $C(\Psi )$ as ${\mathcal P} = {\mathcal A} \cap \overline{\mathcal F}$\,.
\end{definition}

\hspace*{3.2cm} \includegraphics[scale=0.4]{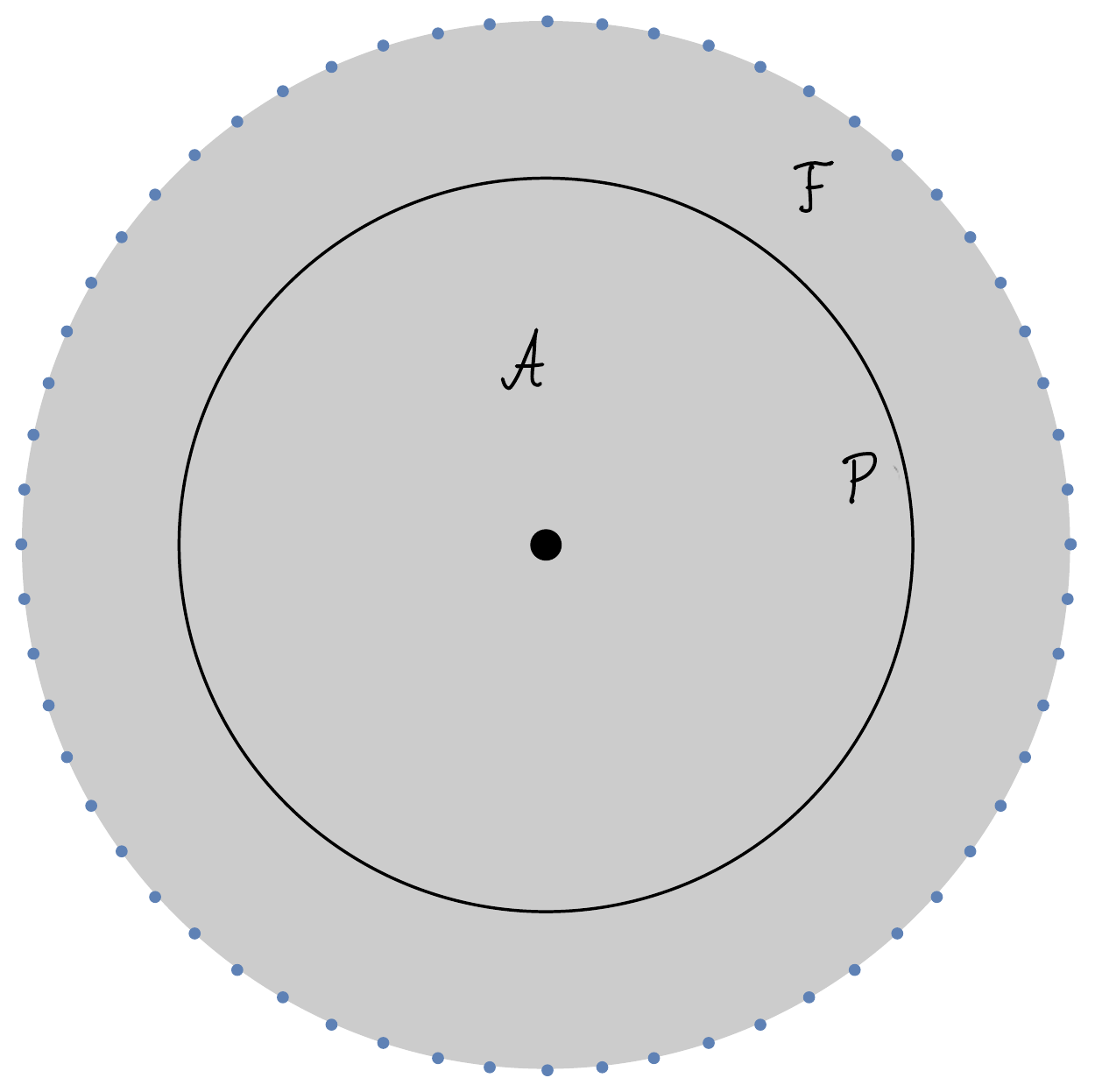} \\

\begin{theorem}\label{theorem-Planck-boundary}
The Planck	 boundary of $C(\Psi)$ is a continuum in $C(\Psi)$  consisting of the minimal ample subcontinua of $\Psi$, which are the terminal pseudo-arcs from the decomposition \eqref{eq:Jones-decomposition}.	 
\end{theorem}

\begin{proof}
By Theorem 19.8 \cite{Illanes-Nadler}, $C(\Psi)$ is unicoherent. Unicoherent means that if we write the continuum as the union of two subcontinua, then their intersection is connected. Therefore, by noticing that ${\mathcal A} \cup \overline{\mathcal F} = C(X)$, we conclude that ${\mathcal P} = {\mathcal A} \cap \overline{\mathcal F}$ is connected.

The Jones decomposition theorem implies that all pseudo-arcs $P_{\alpha}$ from the decomposition \eqref{eq:Jones-decomposition} are ample. However, every proper subcontinuum of a pseudo-arc is terminal and filament in the pseudo-arc. By the terminal property of the pseudo-arcs from the decomposition \eqref{eq:Jones-decomposition},  all proper subcontinua of these pseudo-arcs are terminal and filament in $\Psi$. This shows that the pseudo-arcs from the deccomposition \eqref{eq:Jones-decomposition} are the minimal ample subcontinua and also that  $\mathcal F$ is made of all proper subcontinua of the pseudo-arcs from the decomposition. Therefore, the only boundary points, which can be approximated by elements of $\mathcal F$, are these pseudo-arcs. 
\end{proof}

In the following we use the Whitney map $\mu : 2^{\Psi} \to [0,\mu(\Psi)]$ and  metric $d_{\mu} : 2^{\Psi} \times 2^{\Psi} \to [0, 2 \mu(\Psi)]$, as defined in Section \ref{Section-Whitney}.

Define the Whitney map for $C(\Psi)$ as $\mu_C : C(\Psi) \to [0, \mu (\Psi)]$, which is the restriction of $\mu$ from $2^{\Psi}$ to $C(\Psi)$.

For all $t \in [0,\mu(\Psi)]$ we call $\mu_C^{-1} (t)$ a Whitney level set, and for all
$0 \leq s < t \leq \mu(\Psi)$ we call $\mu_C^{-1} ([s,t])$ a Whitney block.

\smallskip
\smallskip
\smallskip

Based on Proposition \ref{propDecomp} and Theorem \ref{pseudo-circlesinsphere}, we can define a geometric variation of the circle of pseudo-arcs $\Psi$, denoted by $\Psi_0$, which has all minimal ample continua (the maximal pseudo-arcs) on the same Whitney level. 
Indeed, take a continuous decomposition ${\mathcal D}$ of the $2$-sphere $S^2$ into pseudo-arcs with the quotient map $q: S^2 \to S^2$ (see Theorem \ref{pseudo-circlesinsphere}) and a Whitney map $\mu: 2^{S^2}\to [0,\mu(S^2)]$. By Proposition \ref{propDecomp}  we may assume $\mu(D)$ is constant and equals some  fixed positive number $l$ for all $D\in {\mathcal D}$. Let $S^1$ be the unit circle embedded in the $2$-sphere $S^2$. Then $\Psi_0=q^{-1}(S^1)$ is a circle of pseudo-arcs by Theorem \ref{pseudo-circlesinsphere}, and it has all the minimal ample continua, the pseudo-arcs $q^{-1}(y)$, $y\in S^1$, satisfying  $\mu(q^{-1}(y))=l$.

Therefore, we will assume that in the case of $\Psi_0$ we have
\begin{equation}\label{eq:uniquePlancklevel}
\mu_C (P_{\alpha}) = l \; \; \text{for all} \; \alpha \in [0, 2\pi) \, .
\end{equation}
We will call this constant $l$ the {\it Planck constant of $\Psi_0$}.
\smallskip

With assumption \eqref{eq:uniquePlancklevel},  Theorem \ref{theorem-Planck-boundary} has the following corollary.

\begin{cor}
The Planck boundary of $\Psi_0$ is the Whitney level set defined by the Planck constant:
$${\mathcal P} = \mu_C^{-1} (l)\, .$$	
\end{cor}

In the next theorem we list the most important properties of $C(\Psi_0)$.
 
\begin{theorem}\label{theorem-main} Consider $C(\Psi_0)$ endowed with the Whitney metric $d_{\mu}$ defined in Definition \ref{def-Whitney-metric}. Then the following properties hold:\\
(a) The length structure defined by the order-arcs has a nonnegative metric curvature.\\
(b) The filament portion has the decomposition
\begin{equation}\label{eq:filament-decomposition} {\mathcal F} = \bigcup_{\alpha \in [0,2\pi)} \left( C(P_{\alpha}) \setminus \{ P_{\alpha} \} \right) \, ,\end{equation}
where each $C({\mathcal P}_{\alpha})$ is uniquely arcwise-connected and nowhere locally connected.\\
 (c) In the filament portion, two points from different maximal pseudo-arc hyperspaces are connected by concatenations of order arcs going through the interior of the ample region.\\
(d) If $0 < t < l$, then $\mu_C^{-1} (t)$ is a subcontinuum of $C(\Psi_0)$, which is homeomorphic to $\Psi_0$ and $\mu_C^{-1} (t) \subset {\mathcal F}$.\\
(e) The ample region ${\mathcal A}$ is a topological planar disk. In particular, it is contractible, locally contractible and locally connected.\\
(f) If $l \leq t < \mu(\Psi_0)$ then $\mu_C^{-1} (t)$ is a subcontinuum of $C(\Psi_0)$ homeomorphic to $S^1$ with $\mu_C^{-1} (t) \subset {\mathcal A}$.
\end{theorem}

\begin{proof} (a) Using the fact that $\mu_C$ restricted to any order-arc is an isometry, we conclude that if we extend an order-arc to a maximal order-arc connecting a degenerate subcontinuum to $\Psi_0$, the extension will have length equal to $\mu_C (\Psi_0)$. This implies that in any triangle formed by order arcs, the sum of lengths of the two shorter sides equal the length of the third side. Therefore the metric curvature is nonnegative (see Definition 4.1.9 in \cite{Burago}).

(b) The proof of \eqref{eq:filament-decomposition} follows from the facts described in the proof of Theorem \ref{theorem-Planck-boundary}. The terminal property of the pseudo-arcs in the decomposition \eqref{eq:Jones-decomposition} implies that for two disjoint pseudo-arcs, their hyperspaces cannot intersect. As for any hereditarily indecomposable continuum, each $C({\mathcal P}_{\alpha})$ is uniquely arcwise-connected by order arcs. 

(c) Consider $P_{\alpha} \neq P_{\beta}$ from the decomposition \eqref{eq:Jones-decomposition}, $p \in C(P_{\alpha}) \setminus \{ P_{\alpha} \}$  and $q \in C(P_{\beta}) \setminus \{ P_{\beta}$ \}. Then, there are unique order-arcs connecting $p$ to $P_{\alpha}$ and $q$ to $P_{\beta}$. 
As $P_{\alpha}$ and $P_{\beta}$ are on the same Whitney level, the only way connecting them by concatenations of order-arcs is  going through $\Psi_0$ or some other members of of the interior of $C(\Psi_0)$ containing both $P_{\alpha}$ and $P_{\beta}$.

(d) If $0 < t < l$, then by definition, $\mu_C^{-1} (t) \subset \mathcal F$. As shown in Theorem 19.9 \cite{Nadler}, it is also a subcontinuum of $C(\Psi_0)$. Moreover, as hereditarily indecomposable is a Whitney property, for each $0 \leq \alpha < 2\pi$, $\mu_C^{-1} (t) \cap C({\mathcal P}_{\alpha})$ is homeomorphic to ${\mathcal P}_{\alpha}$ and contractible within $C({\mathcal P}_{\alpha})$ to ${\mathcal P}_{\alpha}$. Therefore, $\mu_C^{-1} (t)$ admits the same type decomposition as in Theorem \ref{uniqueness}, which leads to the fact that  $\mu_C^{-1} (t)$ is homeomorphic to $\Psi_0$.

(e) The construction of the circle of pseudo-arcs implies that the ample portion is homeomorphic to the closed unit disk bounded by the Planck boundary.  Therefore, the ample portion ${\mathcal A}$ is contractible, locally contractible and locally connected, as it is also mentioned in Proposition \ref{prop-ample-filament}. 
 
(f) If $l \leq t < \mu (\Psi_0)$, then $\mu_C^{-1} (t) \subset \mathcal A$ and as such, it is homeomorphic to  the Whitney level sets of a circle. We finish the proof by noting that being homeomorphic to a circle is a Whitney property \cite{Nadler}. 

\end{proof}

The following is a schematic picture of the hyperspace of a circle of pseudo-arcs. 

\hspace*{2cm} \includegraphics[scale=0.5]{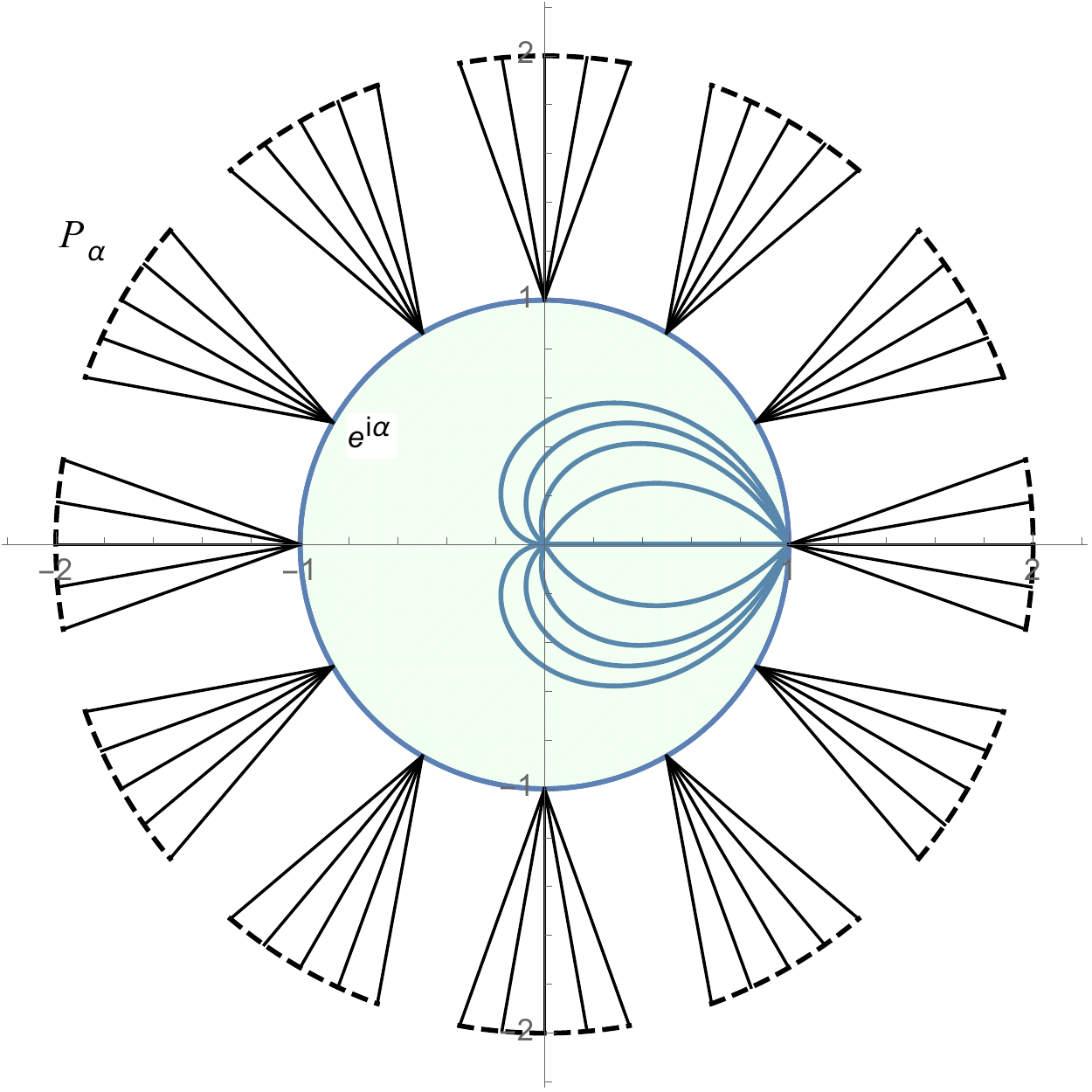} \\

\begin{remark} If we do not assume \eqref{eq:uniquePlancklevel}, not every maximal pseudo-arcs might be on the same level set. Since the Planck boundary is a compact set, $\mu_C$ attains its maximum and minimum on $\mathcal P$. We can define the upper and lower Planck constants of $\Psi$ as
$$L = \max \{ \mu_C (P) \, : \; P \in \mathcal P \} \, ,$$
and
$$l = \min \{ \mu_C (P) \, : \; P \in \mathcal P \} \, .$$	

The connectedness of ${\mathcal P}$ implies that $\mu_C ({\mathcal P}) = [ l, L]$ and
$${\mathcal P} \subset \mu_C^{-1} ([l,L])\, .$$
In this case, the only difference compared to Theorem \ref{theorem-main} is that in $(f)$ we need to use  $L \leq t < \mu(\Psi)$.
\end{remark}




\section{Some philosophical questions and conclusions}

\begin{remark} 
The circle and its rotations are a universal model for measuring time. Can the circle of pseudo-arcs be used as a model for such purpose?

Consider an isotopy (continuous set)  $r_t :  S^1 \to S^1$, $t \in [\alpha, \beta]$,  of rotations of the unit circle $S^1$ by  angles  $t$.  Each $r_t$ can be lifted, with respect to the map $\Pi$ from \eqref{eq:Jones-decomposition}, to a homeomorphism $R_t$ of the circle of pseudo-arcs $\Psi$ \cite{Lewis-85}. This homeomorphism preserves the maximal pseudo-arcs, which leads to a homeomorphisms $\widehat{R}_t:  C(\Psi) \to C(\Psi)$ of the hyperspace $C(\Psi)$.  However, the set of homeomorphisms $\widehat{R}_t$, $t \in [\alpha, \beta]$ displays discontinuities in the filament portion because of the homotopical rigidity of the pseudo-arc as mentioned in Section \ref{Section-Duality}. Therefore, any such attempted lifting of an isotopy of $r_t$, $t \in [\alpha, \beta]$, to an isotopy $\widehat{R}_t$, $t \in [\alpha, \beta]$,  would fail and remain continuous only in the ample portion but discontinuous on the filament portion. Moreover, using the homogeneity of the circle of pseudo-arcs, we can make these discontinuities to be small jumps only. 

So, how would time look like if measured by the circle of pseudo-arcs? On larger scale it would be similar to our usual model of the rotations of the circle. On a small scale it would have to involve ``quantum" jumps. 
Remarkably, our subjective perception of time not always follows the usual smooth uniform model of rotations, but also  involves the sensation of shifts or ``discontinuities."
\end{remark}

\begin{remark}
The structure of our Universe reveals the existence of a sub-Planck portion, then a boundary which shows quantum structure, and finally the macroscopic world. The hyperspace of circles  of pseudo-arcs gives a simple mathematical model of the three portions living together. The filament portion represents the sub-Planck portion, in which there is no local connectedness, all sub-particles are homeomorphic and indistinguishable. As they connect to each-other, they reach a critical size, and become the "elementary-particle" pseudo-arcs in the Planck boundary. These pseudo-arcs behave like point-like objects and can be arcwise-connected by concatenations of arcs (or strings) moving in-an-out of the ample partition. So, from the view-point of the ample region, the Planck boundary can be seen only by the the so-called ``quantum-effect." The ample partition symbolizes the macroscopic world, where more complex assembly of particles is possible.
\end{remark}

\begin{remark} Consider the Menger curve (also known as  the Menger cube or sponge) which was first described by Karl Menger in 1926. The following picture shows the initial three steps of its construction based on dividing cubes into 27 smaller cubes and removing the middle ones.\\
\smallskip
\smallskip

\hspace*{4cm}\includegraphics[scale=0.2]{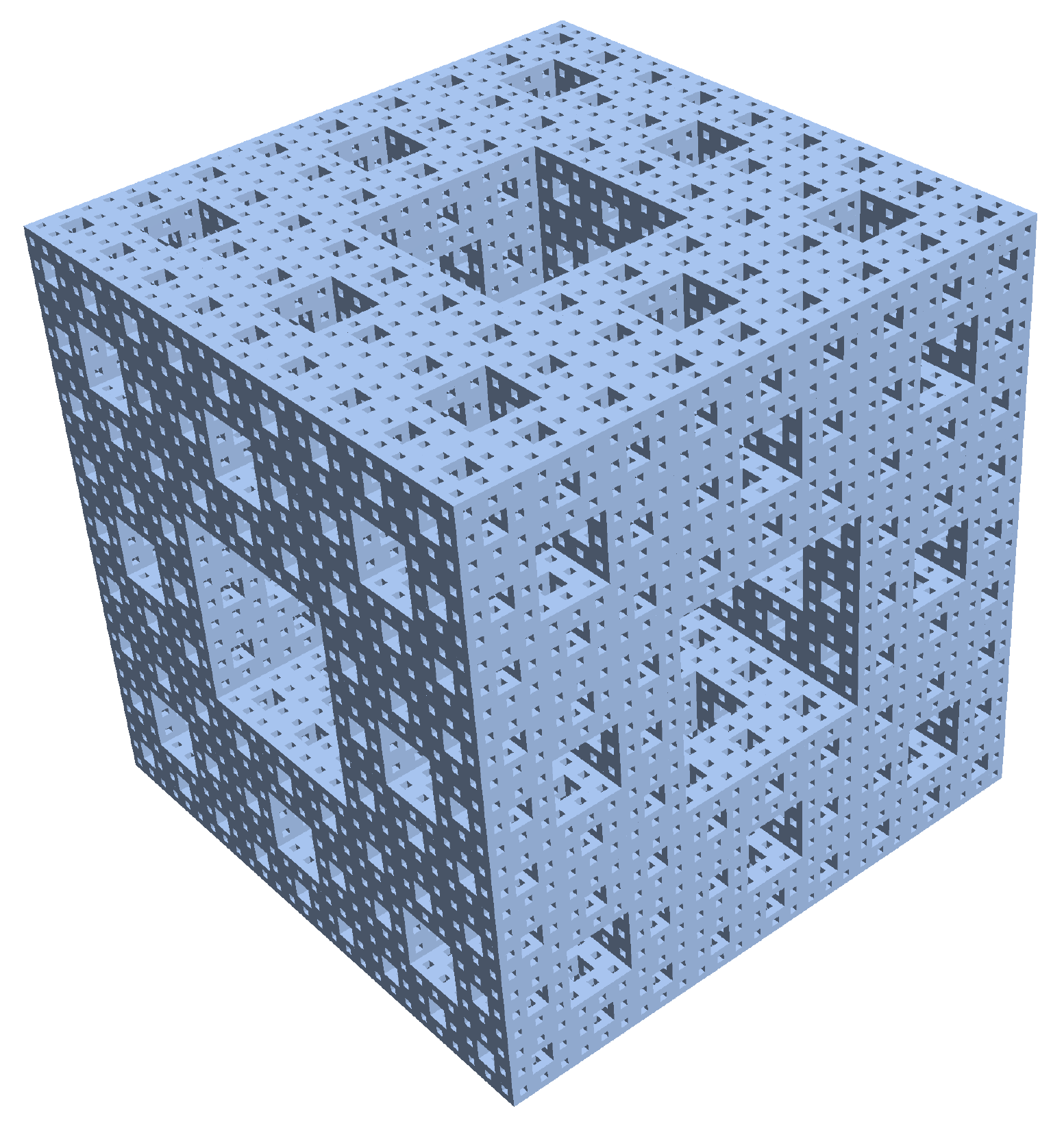}

The Menger curve is a fractal set, with fractal dimension depending on the size of the holes (about 2.7 for the one in the picture), but with topological dimension always equal to one. Using Anderson's topological characterization \cite{Anderson-58}, its geometry can be rearranged to become  any $3$-manifold pierced with countably many  holes of infinitesimally small diameters. The Menger curve is homogeneous. In fact, it has unusually strong homogeneity properties. 

There also exists a {\it Menger curve of pseudo-arcs} \cite{Lewis-85}, a homogeneous continuum admitting a continuous decomposition into terminal pseudo-arcs such that the quotient space is homeomorphic to the Menger curve. An investigation similar to the study of the hyperspace of subcontinua of the circle of pseudo-arcs from the previous section can also be performed, with similar results, on the hyperspace of subcontinua of the Menger curve of pseudo-arcs. In this case, the Planck boundary is a Menger curve, which resembles the ``emptiness'' of space at quantum level. Note that the hyperspace of subcontinua of the Menger curve of pseudo-arcs is infinite dimensional. Can the Menger curve or the Menger curve of pseudo-arcs be used as a model for space? \\

\end{remark}

{\bf Three Centennials.}

\smallskip

About hundred years ago three new and unusual fields of research emerged: quantum mechanics, Jung's synchronicity  and hereditarily indecomposable continua. Seemingly unrelated to each other, they carry remarkable similarities. 

Quantum mechanics is an important component of new physics  developed in the first part of 20th century, which created a change in science on the scale similar to the Newtonian revolution of 17th century. It has important philosophical consequences, too. Quantum mechanics not only proposes the laws of physics at atomic and subatomic scales, but also questions the deterministic nature of the Universe  modeled by classical physics. 

C. G. Jung introduced to Western psychology and philosophy the idea of {\it synchronicity} as dual to the well established principle of {\it causality}. It is easy to draw connection between the determinism of classic physics and causality. Is there any verifiable relation between synchronistic events and quantum phenomena?

The indecomposable and hereditarily indecomposable continua, and in particular the 
pseudo-arc, not only provided  mathematical objects with unexpected properties, but also showed that the majority of connected shapes filling the most common spaces in mathematics are of this form. Moreover, even though the spaces in mathematics are designed to model, analyze  and measure processes in the Universe, the  hereditarily indecomposable continua cannot be measured and analyzed with traditional methods. Can the constructions involving the pseudo-arc, such as the circle of pseudo-arcs and the Menger curve of pseudo-arcs,  be used to model quantum phenomena?
\smallskip
\smallskip
\smallskip



\end{document}